\documentclass{article}
\usepackage{graphicx}

\usepackage[utf8]{inputenc}
\usepackage[T1]{fontenc}
\usepackage[english]{babel}
\usepackage{amsmath,amssymb,amsthm}
\usepackage{graphicx} 
\usepackage[margin=1.in]{geometry}
\usepackage[shortlabels]{enumitem}
\usepackage{todonotes}
\usepackage[bibstyle=numeric, backend=biber, sorting=none, citestyle=numeric-comp, giveninits,url=false,maxbibnames=5,isbn=false,date=year]{biblatex} 
\usepackage{csquotes}\MakeOuterQuote"
\usepackage{soul}\setstcolor{red}
\usepackage[colorlinks,citecolor=blue,linkcolor=blue]{hyperref}
\usepackage[capitalize]{cleveref}
\usepackage{comment}
\usepackage{orcidlink}

\DeclareUnicodeCharacter{0141}{\Lbar{}}

\newtheorem{theorem}{Theorem}
\newtheorem{lemma}[theorem]{Lemma}
\newtheorem{corollary}[theorem]{Corollary}

\theoremstyle{definition}

\newtheorem*{problem*}{Problem}
\newtheorem*{assumption*}{Assumption}

\newtheorem*{warning*}{Warning}

\newcommand{\ip}[2]{\langle #1,#2\rangle}

\newcommand{\ketbra}[2]{|#1\rangle\langle#2|}

\newcommand{\kettbra}[1]{\ketbra{#1}{#1}}

\newcommand{\norm}[1]{\lVert #1\rVert}
\newcommand{\oo}{\infty}
\newcommand{\ox}{\otimes}
\newcommand{\mc}{\mathcal}

\DeclareMathOperator{\tr}{tr}

\renewcommand{\Bar}{\overline}

\newcommand{\hide}[1]{}

\def\B{{\mc B}}

\def\H{{\mc H}}

\def\K{{\mathcal K}}

\def\RR{{\mathbb R}}

\DeclareMathOperator{\lin}{span}

\newcommand{\placeholder}[0]{{\,\cdot\,}}

\newcommand{\verteq}{\scalebox{1.125}{\rotatebox{90}{$\ =\ $}}}

\newcommand{\Lbar}{\L{}}

\title{Revisiting the operator extension of strong subadditivity}
\author{Lauritz van Luijk\,\orcidlink{0000-0003-3153-549X}, Alexander Stottmeister\,\orcidlink{0000-0002-0145-0877}, Henrik Wilming\,\orcidlink{0000-0002-0306-7679}}
\date{\small Leibniz Universit\"at Hannover, Institut f\"ur Theoretische Physik, Appelstraße 2, 30167 Hannover, Germany\\  \today}

\begin{document}

\maketitle
\begin{abstract}
\noindent 
We give a new proof of the operator extension of the strong subadditivity of von Neumann entropy
\begin{equation}\label{eq:operator-ssa-old}
    \rho_A \ox \sigma_{BC}^{-1} \geq \rho_{AB} \ox \sigma_C^{-1}
\end{equation}
by identifying the mathematical structure behind it as Connes' theory of spatial derivatives.
This immediately generalizes the inequality to arbitrary inclusions of von Neumann algebras.
In the case of standard representations, it reduces to the monotonicity of the relative modular operator. 
\end{abstract}

\null

Recently, the operator inequality \eqref{eq:operator-ssa-old}
for bipartite quantum states $\rho_{AB}$ (arbitrary) and $\sigma_{BC}$ with full rank on a finite-dimensional Hilbert space $\H_A\ox\H_B\ox \H_C$ was elegantly proven in \cite{lin_new_2023} (see also \cite{carlen_inequality_2024}). 
The inequality implies strong subadditivity of von Neumann entropy, originally proven in \cite{lieb_proof_1973}, and hence corresponds to an operator strengthening of strong subadditivity.
This curious inequality asks for an explanation in terms of a deeper mathematical structure.

It was observed in \cite[Sec.~3]{lin_new_2023} that \eqref{eq:operator-ssa-old} resembles---but does not follow from---the monotonicity of the relative modular operator, which takes the form $\rho \ox \sigma^{-1}$ for two density matrices $\rho,\sigma$.
In the following, we will see that \eqref{eq:operator-ssa-old} is an immediate consequence of Connes' theory of spatial derivatives on von Neumann algebras \cite{connes_spatial_1980},
which can be understood as a generalization of (relative) modular theory to settings where a given von Neumann algebra $M\subset B(\H)$ is not in \emph{standard form}, i.e., there is no cyclic and separating vector $\Omega\in\H$ for $M$. Indeed, \eqref{eq:operator-ssa-old} is simply the well-known monotonicity of the spatial derivative \cite[Lem.~4.2]{ohya}, applied to the inclusions of matrix factors.
This observation has two consequences:

\begin{enumerate}
    \item There exists an appropriate generalization of \eqref{eq:operator-ssa-old} to the setting of von Neumann algebras. In the standard representation, it reduces to the monotonicity of the relative modular operator \cite[Eq.~(II.1.3)]{borchers_revolutionizing_2000}.
    
    \item We find the following equivalent formulation of \eqref{eq:operator-ssa-old}
      \begin{align}\label{eq:new-ssa}
        \tr_{C}\big(\sigma_C^{-1/2}X_{ABC}\sigma_C^{-1/2}\big)\leq \tr_{BC}\big(\sigma_{BC}^{-1/2} X_{ABC}\sigma_{BC}^{-1/2}\big)
    \end{align}
    as an operator inequality valid for all positive operators $X_{ABC}$ (we suppress explicit identities).
\end{enumerate}

To see the equivalence of \eqref{eq:operator-ssa-old} with \eqref{eq:new-ssa}, we evaluate the latter with $X_{ABC} = \ketbra\xi\xi$ an arbitrary rank-one projection, and take expectation values of the resulting operator on $AB$ against an arbitrary state $\rho_{AB}$:
\begin{eqnarray} 
 \tr\big(\rho_{AB}\,\tr_C(\sigma_C^{-1/2} \kettbra\xi \sigma_C^{-1/2})\big) &\leq &\tr\big(\rho_{A}\,\tr_{BC}(\sigma_{BC}^{-1/2} \kettbra\xi \sigma_{BC}^{-1/2})\big) \nonumber\\
 \verteq \quad\quad&&\quad\quad\verteq \nonumber\\[-22pt] 
  \\[3pt]
\langle \xi,\rho_{AB}\ox \sigma_C^{-1}\xi\rangle  &\leq& \langle \xi,\rho_A \ox \sigma_{BC}^{-1}\xi\rangle.\nonumber
\end{eqnarray}
Since these (in-)equalities hold for all $\xi\in\H$ and states $\rho_{AB}$, the equivalence of \eqref{eq:operator-ssa-old} and \eqref{eq:new-ssa} follows.

\paragraph{Funding.} AS and LvL have been funded by a Stay Inspired Grant of the MWK Lower Saxony (Grant ID: 15-76251-2-Stay-9/22-16583/2022).

\section{Finite-dimensional case}

We first present the finite-dimensional version of the general von Neumann algebraic proof of \eqref{eq:operator-ssa-old}.

Consider a tripartite system with Hilbert space $\H=\H_A\ox\H_B\ox\H_C$.
We fix a full-rank state $\sigma_{BC}$ on $\H_{BC}=\H_B\ox \H_C$ and denote by $\Omega\in \H_{BC}\ox \H_{\Bar{BC}}$ the unnormalized maximally entangled state $\sum_i \eta_i \ox \eta_i$, where $\{\eta_i\}$ constitutes an orthonormal basis of eigenvectors of $\sigma_{BC}$ and $\Bar B$, $\Bar C$ denote copies of $B$, $C$, respectively. 
Then $\Omega_\sigma := \sigma_{BC}^{1/2} \Omega$ is a normalized vector that purifies $\sigma_{BC}$.
Let us define the vector spaces 
\begin{align}
    \K_C := \lin \{ a_{C}\Omega_\sigma : a_{C} \in\B(\H_C)\} \subset  \K_{BC} := \lin\{ a_{BC} \Omega_{\sigma} : a_{BC}\in B(\H_{BC})\}. 
\end{align}
Since $\sigma_{BC}$ has full rank, we have $\K_{BC} = \H_{BC}\ox \H_{\Bar{BC}}$.
For each $\xi \in \H$ we now define a linear operator $ R^{X}(\xi): \K_{X} \to \H$ via
\begin{align}
    R^{X}(\xi)a_{X} \Omega_\sigma = a_{X} \xi, \quad a_{X} \in B(\H_X),
\end{align}
where $X=C$  or $X=BC$.
Evidently, $R^{X}(\xi)a_{X} = a_{X}R^{X}(\xi)$ for $a_{X}\in B(\H_X)$. Moreover, $R^C(\xi)$
is simply the restriction of $R^{BC}(\xi)$ to $\K_C$. 
Denoting by $P$ the orthogonal projection onto $\K_C$ (which commutes with the action of $B(\H_{C})$ on $\K_{BC}$), we thus find
\begin{align}
    \theta^{C}(\xi,\xi) := R^{C}(\xi)R^{C}(\xi)^* = R^{BC}(\xi)PR^{BC}(\xi)^* \leq R^{BC}(\xi)R^{BC}(\xi)^* =: \theta^{BC}(\xi,\xi),
\end{align}
where $\theta^{C}(\xi,\xi)\in B(\H_{C})'=B(\H_{AB})$ and $\theta^{BC}(\xi,\xi)\in\B(\H_{BC})'=B(\H_{A})$.
A useful way to explicitly represent $R^X$ and $\theta^X$ is via tensor-network diagrams (which we draw from right to left):
\begin{align}
    R^{BC}(\xi) = \vcenter{\hbox{\includegraphics{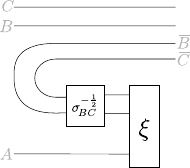}}} \quad\implies \quad \theta^{BC}(\xi,\xi) = \vcenter{\hbox{\includegraphics{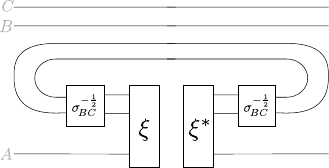}}}\ ,
\end{align}
where an arc denotes an (unnormalized) maximally entangled state.
Inspecting the diagram for $\theta^{BC}(\xi,\xi)$ we see that indeed 
\begin{align}
    \theta^{BC}(\xi,\xi) = \tr_{BC}\big( \sigma_{BC}^{-1/2}\kettbra\xi \sigma_{BC}^{-1/2}\big). 
\end{align}
The isometry $V:\K_{C}\to \K_{BC}$ such that $P = VV^*$ can be written as
\begin{align}
    V = \vcenter{\hbox{\includegraphics{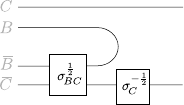}}}\ ,
\end{align}
which leads to 
\begin{align}
    \theta^C(\xi,\xi) = \vcenter{\hbox{\includegraphics{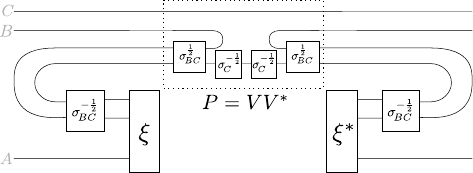}}} = \vcenter{\hbox{\includegraphics{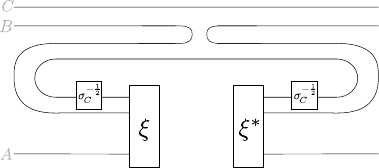}}}\ . \nonumber
\end{align}
Hence
\begin{align}
 \theta^{C}(\xi,\xi) = \tr_{C}\big(\sigma_{C}^{-1/2}\kettbra\xi \sigma_{C}^{-1/2}\big), 
\end{align}
showing \eqref{eq:new-ssa}.
Note  we only needed to use the inequality $APA^* \leq AA^*$, which is true for arbitrary operators $A$ and projections $P$. 

\section{General case}

We now discuss the general case of inclusions of von Neumann algebras $N \subset M$ on a Hilbert space $\H$.
We directly formulate the results in terms of \emph{weights}, which are positive, but generally unbounded, linear functionals on a von Neumann algebra that generalize the notion of a state on a von Neumann algebra. They can be thought of as non-commutative versions of measures on an infinite measure space (e.g., the Lebesgue measure on $\RR$). We will assume familiarity with the basic terminology for weights, see for example \cite{takesaki2}.

Let $\phi,\psi$ be normal semifinite (ns) weights on $N'$ and $M$, respectively.
We assume that the restrictions $\phi|_{M'}$ and $\psi|_{N}$ to the respective subalgebras $N\subset M$ and $M'\subset N'$ are again semifinite (and hence ns weights) and that $\phi$ is faithful.
The setting of the introduction corresponds to $N = B(\H_A) \subset M= B(\H_{A}\ox \H_B) \subset B(\H_A\ox\H_B\ox\H_C)$ and $\psi  = \tr(\rho_{AB}\placeholder)$, $\phi =\tr(\sigma_{BC}\placeholder )$ are normal states on the respective algebras. 

This section aims to explain how the spatial theory of von Neumann algebras, developed by Connes in \cite{connes_spatial_1980}, gives rise to the operator inequality  
\begin{equation}\label{eq:monotonicity-sd}
    \frac{d\psi}{d\phi|_{M'}} \,\leq\, \frac{d\psi|_{N}}{d\phi},
\end{equation}
where the derivatives denote Connes' spatial derivatives, explained below.
In the setting of finite-dimensional quantum information theory discussed above the spatial derivatives are $d\psi/d(\phi|_{M'}) = \rho_{AB}\ox \sigma_C^{-1}$ and $d(\psi|_{N})/d\phi=\rho_A\ox\sigma_{BC}^{-1}$, so that \eqref{eq:monotonicity-sd} reduces to \eqref{eq:operator-ssa-old}.
\\

We need to recall some preliminaries. 
For now, we fix a von Neumann algebra $A$ on a Hilbert space $\H$ together with a normal semifinite faithful (nsf) weight $\omega$.
We will later consider the cases $(A,\omega) = (N',\phi)$ and $(A,\omega)=(M',\phi|_{M'})$, in the notation introduced above.

We denote the GNS representation by $(\pi_\omega,L^2(A,\omega),\eta_\omega)$. 
Here, $L^2(A,\omega)$ is defined as the Hilbert space completion of the left-ideal $\mathfrak{n}_{\omega} = \{ a\in A : \omega(a^*a)<\oo\}\subset A$ with respect to the inner product $(a,b)\mapsto \omega(a^*b)$, $\eta_\omega: \mathfrak{n}_\omega \to L^2(A,\omega)$ denotes the induced injection, and $\pi_\omega$ denotes the representation defined by $\pi_\omega(a)\eta_\omega(b) = \eta_{\omega}(ab)$.
If $\omega$ is a normal faithful state, then $1\in\mathfrak{n_\omega}$, the GNS vector is given as $\Omega_\omega=\eta_\omega(1)$ and $\eta_\omega=\pi_\omega(\placeholder)\Omega_\omega$.

Following \cite{connes_spatial_1980}, we use the GNS representation to define, for each $\xi\in\H$, an operator $R^{\omega}(\xi) : L^2(A,\omega) \supset D(R^{\omega}(\xi)) \to \H$ by 
$D(R^{\omega}(\xi)) =\eta_{\omega}(\mathfrak{n}_{\omega})$, $R^{\omega}(\xi) \eta_{\omega}(a) = a\xi$, $a\in\mathfrak{n}_\omega$.
We note the following intertwining relation:
\begin{equation}\label{eq:intertwining-rel}
    aa'R^\omega(\xi) = R^\omega(a'\xi)\pi_\omega(a), \qquad a\in A,\ a'\in A'.
\end{equation}
A vector $\xi$ is called $\omega$-bounded if $R^{\omega}(\xi)$ is a bounded operator, i.e., there is a  constant $C_{\xi}>0$ such that $\norm{a\xi}\leq C_{\xi}\, \omega(a^*a)^{1/2}$ for all $a\in A$. 
The set of $\omega$-bounded vectors is a dense subspace in $\H$ \cite[Lemma 11.2]{hiaiLecturesSelectedTopics2021}. If $\xi$ is $\omega$-bounded we denote the bounded extension to an operator $L^2(A,\omega) \to \H$ again by $R^{\omega}(\xi)$, and we set
\begin{equation}
    \theta^{\omega}(\xi_1,\xi_2) = R^{\omega}(\xi_1)R^{\omega}(\xi_2)^*.
\end{equation}
for a pair $\xi_1,\xi_2$ of $\omega$-bounded vectors. 
Note that $\theta^\omega$ is linear in the first and conjugate linear in the second entry.
Importantly, \eqref{eq:intertwining-rel} implies that $\theta^\omega(\xi_1,\xi_2)$ is an element of $A'$.
Let us also note that $\theta^{\omega}(\xi,\xi)\ge0$.

Before we describe the spatial derivatives appearing in \eqref{eq:monotonicity-sd}, we prove the von Neumann algebraic version of \eqref{eq:new-ssa} from which we will later infer \eqref{eq:monotonicity-sd}.
In fact, we will see that both inequalities are equivalent.

\begin{lemma}[{\cite[Lem.~4.2]{ohya}}]\label{lem:monotonicity-theta}
    Let $B\subset A$ be a von Neumann subalgebra such that $\omega|_{B}$ is semifinite (this is automatic if $\omega$ is a normal faithful state).
    Then every $\omega$-bounded vector $\xi\in\H$ is $\omega|_{B}$-bounded and we have
    \begin{equation}\label{eq:monotonicity-theta}
        \theta^{\omega|_{\hspace{-.5pt}B}}(\xi,\xi)\le \theta^{\omega}(\xi,\xi).
    \end{equation}
\end{lemma} 

To prove \cref{lem:monotonicity-theta}, we note that in the setting of the Lemma, the inclusion $B\hookrightarrow A$ induces an inclusion on the level of GNS spaces $L^2(B,\omega|_{B})\hookrightarrow L^2(A,\omega)$.
We can regard $L^2(B,\omega|_{B})$ as a closed subspace of $L^2(A,\omega)$.
We have $\eta_{\omega|_{B}}(a)= P\eta_\omega(a)$ for $a\in \mathfrak{n}_{\omega|_{B}}$ and $\pi_{\omega|_{B}}(b) = P\pi_{\omega}(b) = \pi_{\omega}(b)P$, $b\in B$, where $P\in\pi_{\omega}(B)'$ denotes the orthogonal projection of $L^2(A,\omega)$ onto $L^2(B,\omega|_{B})$.
It thus follows that $R^{\omega|_{B}}(\xi) = R^{\omega}(\xi)P$, which shows that $\xi$ is $\omega|_{B}$-bounded if it is $\omega$-bounded.

The monotonicity \eqref{eq:monotonicity-theta} is now an immediate consequence:
\begin{equation*}
    \theta^{\omega|_{\hspace{-.5pt}B}}(\xi,\xi) = R^{\omega|_{\hspace{-.5pt}B}}(\xi)R^{\omega|_{\hspace{-.5pt}B}}(\xi)^* = R^{\omega|_{\hspace{-.5pt}B}}(\xi)  P R^{\omega|_{\hspace{-.5pt}B}}(\xi)^* \le  R^{\omega}(\xi)R^{\omega}(\xi)^* = \theta^\omega(\xi,\xi).
\end{equation*}
This concludes the proof of \cref{lem:monotonicity-theta}.
We note that \cref{eq:monotonicity-theta} remains valid without assuming the faithfulness of $\omega$ (see \cite[Lem.~4.2]{ohya} for further details).

We now come to the definition of Connes' spatial derivatives.
For this, we consider a von Neumann algebra $M$, and we consider an ns weight $\psi$ on $M$ and an nsf weight $\phi$ on $M'$. 
We obtain a quadratic form $(\xi_1,\xi_2)\mapsto \psi(\theta^\phi(\xi_2,\xi_1))$, defined on the subspace of $\phi$-bounded vectors.
This form is closable, and the positive self-adjoint operator corresponding to its closure is the spatial derivative $d\psi/d\phi$:
\begin{align}
    \ip{\xi_1}{\frac{d\psi}{d\phi} \xi_2} = \psi(\theta^\phi(\xi_2,\xi_1)). 
\end{align}
The spatial derivatives are useful tools that allow one to apply modular theory in other representations than the standard representation, in which the spatial derivative reduces to the relative modular operator:
\begin{equation*}
    \frac{d\psi}{d\phi} = \Delta_{\psi | \phi^\top},
\end{equation*}
where $\phi^\top$ denotes the nsf weight $x \mapsto \phi(Jx^*J)$ on $M$, where $J$ is the modular conjugation associated with the standard form of $M$ \cite{haagerup_standard_1975}, which fulfills $JMJ=M'$.

Applying \cref{lem:monotonicity-theta} to the case $A=N'$, $B=M'$, $\omega =\phi$, and evaluating the resulting operator inequality with the weight $\psi$, we get:

\begin{corollary}
    Let $N\subset M$ be von Neumann algebras on a Hilbert space $\H$. 
    Let $\psi$ be an ns weight on $M$ and let $\phi$ be an nsf weight on $N'$.
    If the restrictions $\psi|_{N}$ and $\phi|_{M'}$ are semifinite, then \eqref{eq:monotonicity-sd} holds.
\end{corollary}

We have deduced \eqref{eq:monotonicity-sd} from \eqref{eq:monotonicity-theta}.
As a final remark, we note that the two inequalities are, in fact, equivalent. Indeed, in the setting of \cref{lem:monotonicity-theta}, \cref{eq:monotonicity-sd} implies (cf. \cite[Eq.~(9.1)]{ohya})
\begin{equation}
    \ip\chi{\theta^{\omega|_{\hspace{-.5pt}B}}(\xi,\xi)\chi} = \ip\xi{\frac{d\ip\chi{\placeholder\chi}|_{B'}}{d\omega|_{B}} \xi} \le \ip\xi{\frac{d\ip\chi{\placeholder\chi}|_{A'}}{d\omega} \xi} = \ip\chi{\theta^\omega(\xi,\xi)\chi}, \qquad \chi \in \H,
\end{equation}
and, thus, \eqref{eq:monotonicity-theta} implies \eqref{eq:monotonicity-sd} as well.

\printbibliography

\end{document}